\documentclass[12 pt]{amsart}

\usepackage{amssymb}
\usepackage{amsbsy}
\usepackage{amscd}
\usepackage{amsmath}
\usepackage{amsthm}

\newtheorem{theorem}{Theorem}[section]
\newtheorem{defn}[theorem]{Definition}
\newtheorem{lem}[theorem]{Lemma}
\newtheorem{prop}[theorem]{Proposition}
\newtheorem{rem}[theorem]{Remark}
\newtheorem{example}[theorem]{Example}
\newtheorem{corollary}[theorem]{Corollary}
\newtheorem{pbm}[theorem]{Problem}

\newcommand{\ad}{\operatorname{ad}}
\newcommand{\Cat}{\operatorname{Cat}}
\newcommand{\Coker}{\operatorname{Coker}}
\newcommand{\Der}{\operatorname{Der}}

\newcommand{\Ker}{\operatorname{Ker}}
\newcommand{\Mor}{\operatorname{Mor}}
\newcommand{\mult}{\operatorname{mult}}

\newcommand{\pro}{\operatorname{pr}}
\newcommand{\Proj}{\operatorname{Proj}}
\newcommand{\Spec}{\operatorname{Spec}}

\newcommand{\A}{\mathcal{A}}

\newcommand{\calE}{\mathcal{E}}

\newcommand{\calH}{\mathcal{H}}

\newcommand{\calO}{\mathcal{O}}
\newcommand{\calT}{\mathcal{T}}

\newcommand{\bC}{\mathbb{C}}
\newcommand{\bF}{\mathbb{F}}
\newcommand{\bK}{\mathbb{K}}

\newcommand{\bP}{\mathbb{P}}

\newcommand{\bR}{\mathbb{R}}
\newcommand{\bZ}{\mathbb{Z}}

\newcommand{\frakh}{\mathfrak{h}}
\newcommand{\frakg}{\mathfrak{g}}
\newcommand{\frakS}{\mathfrak{S}}

\newcommand{\m}{{\bf m}}

\title[Arrangements and multiderivations]
{Arrangements, multiderivations, and adjoint quotient 
map of type ADE}

\author[M. Yoshinaga]{Masahiko Yoshinaga}

\address{Department of Mathematics, Kyoto University, Kyoto 606-8502, Japan}

\email{mhyo@math.kyoto-u.ac.jp}

\begin{document}

\begin{abstract}
The first part of this paper is a 
survey on algebro-geometric aspects of sheaves of 
logarithmic vector fields of hyperplane arrangements. 
In the second part we prove that the 
relative de Rham cohomology (of degree two) of 
ADE-type adjoint quotient map is 
naturally isomorphic to the module of certain 
multiderivations. The isomorphism is obtained 
by the Gauss-Manin derivative of the Kostant-Kirillov 
form. 
\end{abstract}

\maketitle

\section{Introduction}

We begin with an example to 
illustrate how the structure of the 
module of logarithmic vector fields $D(\A)$ is 
related to combinatorial problems of a hyperplane 
arrangement $\A$. 
Let $\A$ be a collection 
$\{H_{ij}\mid 1\leq i < j\leq n\}$ of 
hyperplanes $H_{ij}=\{(x_1, \dots, x_n)\in\bK^n
\mid x_i-x_j=0\}\subset \bK^n$, 
where $\bK$ is a fixed field. 
According to the field $\bK$, several enumerative problems 
appear 
for the complement $M(\A)=\bK^n\setminus\bigcup H_{ij}$. 
\begin{itemize}
\item[(i)] If $\bK=\bF_q$ is a finite filed, then the 
complement $M(\A)$ is a finite set, of 
$|M(\A)|=q(q-1)(q-2)\dots(q-n+1).$ 
\item[(ii)] If $\bK=\bR$ is the real numbers, then each 
connected component of $M(\A)$ (the chamber) is 
expressed by the inequality 
$x_{i_1}<x_{i_2}<\dots<x_{i_n}$, where 
$(i_1, \dots, i_n)$ is a permutation of $(1, \dots, n)$. 
There are $n!$ chambers. 
\item[(iii)] If $\bK=\bC$ is the complex numbers, then 
$M(\A)$ is an affine complex variety of $\dim_\bC=n$. 
Using the fibration $(x_1, \dots, x_n)\longmapsto 
(x_1, \dots, x_{n-1})$ and the Leray-Hirsch theorem, 
the Poincar\'e polynomial is computed as 
$\sum_i b_i(M(\A))t^i=(1+t)(1+2t)\dots(1+(n-1)t)$.
\end{itemize}
The formulas in (i)-(iii) are similar in appearance. 
The general theory of arrangements 
\cite{ot} tells us that 
these invariants are combinatorial. 
Namely, they are determined from the poset 
$L(\A)$ of subspaces obtained as intersections. 
Computations of the characteristic polynomial 
$\chi(\A, t)\in\bZ[t]$ 
unify these enumerative problems.

We also consider derivations 
$$
\delta_p=\sum_{i=1}^n x_i^p{\partial_i}, 
$$
($\partial_i=\frac{\partial}{\partial x_i}$) with 
$p=0, 1, \dots, n-1$. 
These satisfy 
\begin{equation}
\label{eq:sai1}
\delta_p(x_i-x_j)=x_i^p-x_j^p\in 
(x_i-x_j), 
\end{equation}
for all $i, j$ and the determinant of coefficients 
\begin{equation}
\label{eq:sai2}
\det
\left(
\begin{array}{ccccc}
1&x_1&x_1^2&\dots&x_1^{n-1}\\
1&x_2&x_2^2&\dots&x_2^{n-1}\\
\vdots&\vdots&\vdots&\ddots&\vdots\\
1&x_n&x_n^2&\dots&x_n^{n-1}
\end{array}
\right)=
\prod_{1\leq i<j\leq n}(x_j-x_i)
\end{equation}
is the product of defining equations. 
The properties (\ref{eq:sai1}) and 
(\ref{eq:sai2}) guarantee that the module 
$$
D(\A)=
\{\delta\in\Der_S\mid \delta(x_i-x_j)\in 
(x_i-x_j), \forall i,j\}
$$
is a free module over $S=\bK[x_1, \dots, x_n]$ 
with basis $\delta_0, \dots, \delta_{n-1}$ 
(this is Saito's criterion \cite{sai-log}). 
Remarkably, the decomposition of $D(\A)$ into a 
direct sum of rank one free modules implies 
the product formulas (i)--(iii) above 
(Terao's factorization theorem \cite{ter-fact}). 
More generally, the algebraic structure of $D(\A)$ 
determines the characteristic polynomial $\chi(\A, t)$ 
by Solomon-Terao's formula 
\cite{st-stf} (see also \S\ref{subsec:log} below). 

The graded $S$-module $D(\A)$ can also be considered 
as a coherent sheaf $\widetilde{D}(\A)$ on 
projective space  $\bP^{n-1}$. 
This fact enables us to employ algebro-geometric 
methods to study $\A$. The structures of these sheaves 
contain information on $\A$. 

The purpose of this paper is to survey algebro-geometric 
aspects of $D(\A)$ and give some related results. 
The paper is organized as follows. 
In \S\ref{sec:arrag}, we start with recalling 
basic notions on logarithmic vector fields for a Cartier 
divisor. We also introduce the module $D(\A, \m)$ of 
logarithmic vector fields 
for an arrangement with multiplicity (multiarrangements) 
in \S\ref{subsec:log}. 
In general, the logarithmic vector fields for multiarrangements 
are much more difficult to analyze than simple arrangements. 
However, freeness of rank $\ell$ simple arrangements is 
closely related to freeness of rank $\ell-1$ multiarrangements. 
We will describe 
freeness criteria for these objects in 
\S\ref{subsec:char1}--\S\ref{subsec:char2}. 
In \S\ref{subsec:curv}, we give a new necessary condition 
for a $3$-dimensional arrangement to be free, in terms of 
plane curves. 
In \S\ref{sec:coxeter}, we will review results on 
freeness of Coxeter multiarrangements. Coxeter arrangements 
are the 
best understood class of multiarrangements. 
In \S\ref{sec:appl}, we will give two applications 
of freeness of Coxeter multiarrangements. 
The first concerns the adjoint quotient map 
$\chi:\frakg\rightarrow\frakg//\ad(G)$ of a simple 
Lie algebra $\frakg$ of $ADE$-type. To describe the 
relative de Rham cohomology of $\chi$, 
the module $D(\A, \m)$ naturally appears. 
In the second application, we will give another 
proof for the freeness of $A_n$-Catalan arrangements, 
which was first proved by Edelman and Reiner \cite{ede-rei}. 

%In \S\ref{subsec:adj}, 
%In \S\ref{subsec:cat}, 

{\em Acknowledgements.} 
The author deeply thanks Professor Kyoji Saito. 
Parts of this article (especially 
\S\ref{subsec:adj}) were done under his supervision. 
This work was supported by 
JSPS Grant-in-Aid for Young Scientists (B) 
20740038.

\section{Algebraic Geometry of logarithmic vector fields}
\label{sec:arrag}

\subsection{Sheaf of logarithmic vector fields}
\label{subsec:sheaf}

Let $X$ be a smooth complex variety and $D\subset X$ 
a Cartier divisor. 
Let $U\subset X$ be an open subset of $X$. Suppose that 
there exists $h\in\Gamma(U, \calO_X)$ such that 
$U\cap D=\{h=0\}$. 
Let $\delta\in\Gamma(U, \calT_X)$ be a 
section of the tangent sheaf on an open subset 
$U\subset X$ (i.e., a holomorphic vector 
field on $X$). 
The section $\delta$ is said to be 
{\em logarithmic tangent to $D$} if 
$\delta h\in h\cdot \calO_U$. 
The sheaf of vector fields logarithmic tangent to $D$ 
is denoted by $\calT_X(-\log D)$. 
The sheaves of logarithmic forms are also similarly 
defined as 
$$
\Omega^p_X(\log D)=\{
\omega\in\frac{1}{h}\Omega^p_X\mid 
\omega\wedge dh\mbox{ is holomorphic }\}. 
$$
They were introduced by K. Saito in \cite{sai-log}. 
He proved that 
they are reflexive sheaves and if $\calT_X(-\log D)$ 
(or $\Omega^1_X(\log D)$) is locally free 
then $\Omega^p_X(\log D)=\bigwedge^p \Omega^1_X(\log D)$. 
We also note that if $\dim X=2$, $\calT_X(-\log D)$ 
is a locally free sheaf. 

\begin{example}
Let $X=\bP_\bC^2=\Proj\bC[z_0, z_1, z_2]$. 
Using the Euler sequence (\cite{oss-vb}) 
$$
0
\longrightarrow
\calO_{\bP^2}
\longrightarrow
\calO_{\bP^2}(1)^3
\longrightarrow
\calT_{\bP^2}
\longrightarrow
0, 
$$
we have the following. 
\begin{itemize}
\item[(1)] If $D_0=\{z_0=0\}\subset\bP^2$. Then 
$\calT_{\bP^2}(-\log D_0)\cong\calO_{\bP^2}(1)^2$. 
\item[(2)] If $D_1=\{z_0z_1=0\}\subset\bP^2$. Then 
$\calT_{\bP^2}(-\log D_1)\cong\calO_{\bP^2}(1)\oplus\calO_{\bP^2}$. 
\item[(3)] If $D_2=\{z_0z_1z_2=0\}\subset\bP^2$. Then 
$\calT_{\bP^2}(-\log D_2)\cong\calO_{\bP^2}^2$. 
\item[(4)] If $D_3=\{z_0^2+z_1^2+z_2^2=0\}\subset\bP^2$. Then 
$\calT_{\bP^2}(-\log D_3)\cong\calT_{\bP^2}(-1)$. 
(Sketch: $\bigoplus_d\Gamma(\calT_{\bP^2}(-\log D_3)(d))$ 
is generated by 
$
\delta_0=z_2\partial_1-z_1\partial_2, 
\delta_1=-z_2\partial_0+z_0\partial_2$ and 
$\delta_2=z_1\partial_0-z_0\partial_1$ with a relation 
$z_0\delta_0+z_1\delta_1+z_2\delta_2=0$ this induces a 
resolution which is isomorphic to shifted Euler sequence.) 
\end{itemize}
\end{example}
The examples above, $\calT_X(-\log D)$ is always 
a uniform sheaf. However for ``generic'' divisors of higher 
degrees, we obtain ``generic'' sheaves. We can sometimes 
recover the original divisor $D$ from the sheaf 
$\calT_{\bP^n}(-\log D)$. 
(Dolgachev and Kapranov called this type of 
result a ``Torelli-type'' theorem in \cite{dk}.) 
Let us recall two results in this direction. 
First one is due to Dolgachev and Kapranov, concerning 
the case of a union of generic hyperplanes.

\begin{theorem}
\cite{dk} 
Let $m\geq 2n+3$ and $\A_i=\{H_{i1}, H_{i2}, \dots, 
H_{im}\}$ ($i=1, 2$) be arrangements of 
generic $m$ hyperplanes $H_{ik}\subset \bP_\bC^n$ in 
$n$-dimensional projective space. 
We denote the union by $\cup\A=\bigcup_{H\in\A}H$. 
If 
$$
\calT_{\bP^n}(-\log(\cup\A_1))\cong
\calT_{\bP^n}(-\log(\cup\A_2)), 
$$
then 
$\cup\A_1=\cup\A_2$. 
\end{theorem}
For smooth divisors $D\subset\bP^n$ defined 
by a homogeneous polynomial $\{f(z_0, \dots, z_n)=0\}$, 
Torelli-type results are related to the following 
property. 
\begin{defn}
The homogeneous polynomial $f(z_0, \dots, z_n)$ 
is said to be {\em Thom-Sebastiani type} if 
there exists a linear transformation 
$A:\bC^{n+1}\rightarrow\bC^{n+1}$ such that 
$f(A(z_0, \dots, z_n))=g(z_0, \dots, z_k)+
h(z_{k+1}, \dots, z_n)$ 
for some $0\leq k\leq n-1$. 
\end{defn}

\begin{theorem}
\cite{uy-cub, uy-sm}
\begin{itemize}
\item[(i)] Let $D_1, D_2\subset\bP^n$ be 
degree $d$ smooth divisors which are 
not Thom-Sebastiani type. Then 
$\calT_{\bP^n}(-\log D_1)\cong
\calT_{\bP^n}(-\log D_2)$ if and only if $D_1=D_2$. 
\item[(ii)] Let $D_1, D_2\subset\bP^2$ be smooth cubic 
curves with non-zero $j$-invariant $j(D_i)\neq 0$. Then 
$\calT_{\bP^n}(-\log D_1)\cong
\calT_{\bP^n}(-\log D_2)$ if and only if $D_1=D_2$. 
\end{itemize}
\end{theorem}

\subsection{Log vector fields for multiarrangements}
\label{subsec:log}

The main theme of this paper is logarithmic vector fields 
for arrangements of hyperplanes. 
Freeness is one of the important properties for 
arrangements. Ziegler \cite{zie-multi} 
showed that a free arrangement 
of rank $\ell$ induces several free multiarrangements of 
rank $\ell-1$ (see \S\ref{subsec:char2}). This means that 
freeness of multiarrangements will be necessary for 
that of simple arrangements. Recently, 
several results on free simple arrangements have been generalized 
to free multiarrangements.

Let $V=\bC^\ell$ be a complex vector space with 
coordinate $(x_1, \cdots, x_\ell)$, $\A=\{H_1, \ldots, H_n\}$ 
be an arrangement of hyperplanes. Let us denote by 
$S=\bC[x_1, \ldots,x_\ell]$ the polynomial ring and fix 
$\alpha_i\in V^*$ a defining equation of $H_i$, i.e., 
$H_i=\alpha_i^{-1}(0)$. We also put 
$Q(\A, \m)=\prod_{i=1}^n\alpha_i^{\m(H_i)}$ and $|\m|=\sum_i \m(H_i)$. 

\begin{defn}
A multiarrangement is a pair $(\A, \m)$ of an arrangement 
$\A$ with a map $ \m: \A\rightarrow \bZ_{\geq 0}$, called 
the multiplicity.  
\end{defn}
An arrangement $\A$ can be identified with a multiarrangement 
with constant multiplicity $m\equiv 1$, which is sometimes called 
a {\em simple arrangement}. 
With this notation, the main object is 
the following module of $S$-derivations which has contact 
to each hyperplane of order $m$. 
\begin{defn}
Let $(\A, \m)$ be a multiarrangement, and define 
$$
D(\A, \m)=\{\delta\in\Der_S| \delta\alpha_i\in(\alpha_i)^{\m(H_i)},
\forall i\}, 
$$
and 
$$
\Omega^p(\A, \m)=\left\{\left.
\omega\in\frac{1}{Q}\Omega^p_V\right| d\alpha_i\wedge\omega
\mbox{ \normalfont{does not have pole along} }H_i, 
\forall i\right\}, 
$$
\end{defn}
The module $D(\A, \m)$ is obviously a graded $S$-module. 
A multiarrangement $(\A, \m)$ is said to be {\em free} with 
exponents $(e_1, \ldots, e_\ell)$ if and only if 
$D(\A, \m)$ is an $S$-free module and there exists a 
basis $\delta_1, \ldots,\delta_\ell\in D(\A, \m)$ such that 
$\det \delta_i=e_i$. Here note that the degree 
$\deg\delta$ of a derivation $\delta$ is the polynomial 
degree, in other words, $\deg (\delta f)=\deg\delta + 
\deg f-1$ for a homogeneous polynomial $f$. 
An arrangement $\A$ is said to be free if $(\A, 1)$ is free. 

Let $\delta_1, \ldots, \delta_\ell\in D(\A, \m)$. Then 
$\delta_1, \cdots, \delta_\ell$ form a $S$-basis of $D(\A, \m)$ 
if and only if 
$$
\delta_1\wedge\cdots\wedge\delta_\ell=
c\cdot Q(\A, \m)\cdot
\frac{\partial}{\partial z_1}\wedge\cdots\wedge\frac{\partial}{\partial z_\ell}, 
$$
where $c\in\bC^*$ is a non-zero constant 
(Saito's criterion \cite{sai-log}). 
From Saito's criterion, we also obtain that if 
a multiarrangement $(\A, \m)$ is free with 
exponents $(e_1, \ldots, e_\ell)$, then 
$|\m|=\sum_{i=1}^\ell e_i$. 

The Euler vector field 
$\theta_E=\sum_{i=1}^\ell x_i\partial_i$ 
is always contained in $D(\A)$. Thus it is 
natural to define 
$D_0(\A):=D(\A)/S\cdot\theta_E$. 
Since $D_0(\A)$ is a graded $S$-module, it 
determines a coherent sheaf 
$\widetilde{D_0}(\A)$ on $\bP^{\ell-1}$. 
On the other hand, 
an arrangement $\A$ defines a Cartier 
divisor $\cup\overline{\A}=\bigcup \overline{H}
\subset\bP^{\ell-1}$. The logarithmic sheaf 
$\calT_{\bP^{\ell-1}}(-\log(\cup\overline{\A}))$ 
determined by the divisor $(\cup\overline{\A})$ 
is related to 
$\widetilde{D_0}(\A)$ by the following formula. 
$$
\calT_{\bP^{\ell-1}}(-\log(\cup\overline{\A}))
\cong
\widetilde{D_0}(\A)[+1]. 
$$
From the sheaf $\widetilde{D}(\A)$, 
we can reconstruct the graded module $D(\A)$ as 
global sections $\Gamma_*(\bP^{\ell-1}, \widetilde{D}(\A)):=
\bigoplus_{d\in\bZ}\Gamma(\bP^{\ell-1}, \widetilde{D}(\A))$. 
More generally, we have the following. 
%See \cite[Lemma ???????]{ay-split} for 
%the proof. 
\begin{prop}
\label{prop:section}
Let $(\A, \m)$ be a multiarrangement. Then 
the natural map 
$$
\Omega^p(\A, \m)\longrightarrow
\Gamma_*(\bP^{\ell-1}, \widetilde{\Omega}^p(\A, \m))
$$
is an isomorphism. 
\end{prop}
\proof
We prove the surjectivity. Since 
$\bigcup_{i=1}^\ell U_i=\bP^{\ell-1}$, where 
$U_i=\{z_i\neq 0\}\subset\bP^{\ell-1}$, is an 
affine open covering, 
any element of right hand side 
$\frac{\omega}{Q}\in
\Gamma_*(\bP^{\ell-1}, \widetilde{\Omega}^p(\A, \m))$ 
can be expressed as 
$$
\frac{\omega}{Q}=
\frac{\omega_1}{z_1^{d_1}Q}=
\frac{\omega_2}{z_2^{d_2}Q}=\dots=
\frac{\omega_\ell}{z_\ell^{d_\ell}Q}, 
$$
with $\frac{\omega_i}{Q}\in\Omega^p(\A, \m)$. 
Using the fact that $S$ is a UFD, it is easily seen 
that $\omega$ is a regular differential form. 
Assume that $z_i$ and $\alpha_H$ are linearly 
independent. Taking the wedge with $d\alpha_H$, 
$d\alpha_H\wedge\frac{\omega_i}{Q}$ does not have 
pole along $H$. So 
$d\alpha_H\wedge\frac{\omega}{Q}=
d\alpha_H\wedge\frac{\omega_i}{z_i^{d_i}Q}$. 
Hence 
$d\alpha_H\wedge\frac{\omega}{Q}\in
\Omega^p(\A, \m)$. 
\qed

Combining the above proposition with a sheaf 
theoretic property of reflexive sheaves, we can 
prove that $\Omega^p(\A, \m)$ is determined by 
$\Omega^1(\A, \m)$ in general. 

\begin{prop}
Assume that 
$\Omega^1(\A_1, m_1)\cong
\Omega^1(\A_2, m_2)$. Then 
$\Omega^p(\A_1, m_1)\cong
\Omega^p(\A_2, m_2)$. 
\end{prop}

\begin{proof}
Let us denote $\calE^p=\widetilde{\Omega}^p(\A, \m)$. 
Let $U:=\bP^{\ell-1}\setminus\bigcup_{H, H'\in\A, H\neq H'}
(\overline{H}\cap \overline{H'})$ 
be the complement to the union of codimension 
$\geq 2$ strata. 
We denote the inclusion map 
$i:U\hookrightarrow\bP^{\ell-1}$. The restriction 
$i^*\calE^p$ is locally free. Hence we have 
$i^*\calE^p=\bigwedge^p i^*\calE^1$. Since $\calE^p$ 
is reflexive, hence normal, we have 
$\calE^p=i_*(i^*(\calE^p))$. Thus 
$\calE^p=i_*(\wedge^p i^*(\calE^1))$. 
(See \cite[\S 1]{har-ref} for basic properties of 
reflexive sheaves.) 
Thus $\calE^p$ is determined by $\calE^1$. Then by 
Proposition \ref{prop:section}, we obtain 
the graded module $\Omega^p(\A, \m)$ from 
its sheafification. 
\end{proof}
The following result shows that $D(\A)$ determines 
the characteristic polynomial $\chi(\A, t)$. 
Corollary \ref{cor:tfact} 
is known as Terao's factorization theorem. 

\begin{theorem}
\label{thm:stf}
\cite{st-stf} 
Denote by $H(\Omega^p(\A), x)\in\bZ[[x]][x^{-1}]$ the 
Hilbert series of the graded module $\Omega^p(\A)$. Define 
\begin{equation}
\label{eq:phi}
\Phi(\A; x, y)=
\sum_{p=0}^\ell H(\Omega^p(\A), x)y^p. 
\end{equation}
Then 
\begin{equation}
\label{eq:stf}
\chi(\A, t)=\lim_{x\rightarrow 1}
\Phi(\A; x, t(1-x)-1). 
\end{equation}
\end{theorem}

\begin{corollary}
\label{cor:tfact}
\cite{ter-fact} 
Suppose that $\A$ is a free arrangement with 
exponents $(e_1, \dots, e_\ell)$. Then 
\begin{equation}
\label{eq:fact}
\chi(\A, t)=\prod_{i=1}^\ell(t-e_i).
\end{equation}
\end{corollary}
The notion of freeness has a geometric interpretation. 
It is equivalent to a splitting $\widetilde{D}(\A)=
\bigoplus_{i=1}^\ell\calO(-e_i)$ of the sheaf 
$\widetilde{D}(\A)$. Then the formula (\ref{eq:fact}) 
indicates that the characteristic polynomial 
is related to the Chern polynomial $c_t(\calE)=
c_0(\calE)+c_1(\calE)t+\dots$ (recall that 
$c_t(\calO(-e))=1-et$) \cite{sil-poi}. 
Indeed, for locally free arrangements, 
Musta\c t\v a and Schenck gave 
a beautiful formula connecting $\chi(\A, t)$ and 
the Chern polynomial. 

\begin{theorem}
\label{thm:ms}
\cite{mus-sch} 
If $\widetilde{D}(\A)$ is a locally free sheaf 
on $\bP^{\ell-1}$, then 
\begin{equation}
\label{eq:ms}
c_t(\widetilde{D_0}(\A))=
t^{\ell-1}\chi_0(\A, 1/t),  
\end{equation}
where $\chi_0(\A, t)=\chi(\A, t)/(t-1)$. 
\end{theorem}
Note that in the case $\ell\leq 3$, the local freeness is 
always satisfied. Thus the Chern polynomial is essentially 
equivalent to $\chi(\A, t)$ and combinatorially computable 
\cite{sch-rk2}.

\subsection{Characterizing freeness}
\label{subsec:char1}

A vector bundle on $\bP^1$ is always a direct 
sum of line bundles (Grothendieck). 
The splitting of vector bundles on $\bP^n (n\geq 2)$ 
is also a well studied subject, e.g., see \cite{oss-vb}. 
There are several criterion to be split. The next 
result is known as Horrocks' criterion. 

\begin{theorem}
\label{thm:hor}
Let $E$ be a rank $r$ holomorphic vector bundle on 
$\bP^n (n\geq 2)$. The following conditions are 
equivalent. 
\begin{itemize}
\item[(i)] $E=\bigoplus_{i=1}^r\calO_{\bP^n}(d_i)$ for some 
$d_1, \dots, d_r\in\bZ$. 
\item[(ii)] $H^i(\bP^n, E(d))=0$ for $\forall 1\leq i\leq n-1$ 
and $\forall d\in\bZ$. 
\item[(iii)] (If $n\geq 3$) $\exists H\subset\bP^n$ a hyperplane 
such that the restriction splits as $E|_H
=\bigoplus_{i=1}^r\calO_{H}(d_i)$. 
\end{itemize}
\end{theorem}

Yuzvinsky \cite{yuz-sheaf, yuz-obst, yuz-loc} developed 
sheaf theory on the intersection lattice $L(\A)$ and 
gave a cohomological criterion for an arrangement $\A$ 
to be free which is similar to Theorem \ref{thm:hor} (ii). 
As an application 
he proved that the set of free arrangements form a 
Zariski open subset in the moduli space of all 
arrangements having the fixed combinatorial type. 

Here we describe a criterion similar to Theorem \ref{thm:hor} (iii). 
We begin with recalling Ziegler's restriction 
\cite{zie-multi}.

Choose a hyperplane $H\in\A$ and coordinate $(z_1, 
\ldots, z_\ell)$ such that $H=\{z_\ell=0\}$. 
Define a submodule $D_0^H(\A)$ of $D(\A)$ 
as follows: 
$$
D_0^H(\A):=
\{\delta\in D(\A)\mid \delta z_\ell=0\}. 
$$

\begin{lem}
\label{lem:zie}
$D(\A)=S\cdot \theta_E\oplus D_0^H(\A)$. 
\end{lem}
\begin{proof}
Let $\delta\in D(\A)$. The assertion is obvious from 
$
\delta=
\left(
\frac{\delta z_\ell}{z_\ell}
\right)\theta_E
+
\left(
\delta-\frac{\delta z_\ell}{z_\ell}\theta_E
\right). 
$
\end{proof}

The arrangement $\A$ determines the restricted 
arrangement $\A^H=\{H\cap H'\mid H'\in\A, H'\neq H\}$ 
on $H$. The restricted arrangement $\A^H$ possesses a 
natural multiplicity 
$$
\begin{array}{cccl}
m^H:&\A^H&\longrightarrow&\bZ\\
&X&\longmapsto&\sharp\{H'\in\A\mid X=H\cap H'\}. 
\end{array}
$$
Ziegler \cite{zie-multi} proved that 
the freeness of $\A$ implies that of $(\A^H, m^H)$. 

\begin{theorem}
\label{thm:ziegler}
\cite{zie-multi}
\begin{itemize}
\item[(1)] 
If $\delta\in D_0^H(\A)$, then $\delta|_{z_\ell=0}\in D(\A^H,
m^H)$. 
\item[(2)] 
If $\A$ is free with exponents $(1, e_2, \dots, e_\ell)$, 
then $(\A^H, m^H)$ is free with 
exponents $(e_2, \ldots, e_\ell)$. 
\end{itemize}
\end{theorem}
\begin{corollary}
\label{cor:zie}
$\A$ is free with exponents $(1, e_2, \dots, e_\ell)$ 
if and only if the following are satisfied. 
\begin{itemize}
\item $(\A^H, m^H)$ is free with exponents $(e_2, \ldots, e_\ell)$. 
\item The restriction induces the surjection 
$D_0^H(\A)\longrightarrow D(\A^H,m^H)$. 
\end{itemize}
\end{corollary}
Using Corollary \ref{cor:zie}, we can establish a 
Horrocks' type criterion for freeness. 
Namely, we will characterize 
freeness by using the freeness of the 
restriction $D(\A^H,m^H)$. 
We first consider the case $\ell=3$. 
By analyzing the Hilbert series of these graded modules 
using the restriction map and Solomon-Terao's formula 
(Theorem \ref{thm:stf}), 
we have the following. 

\begin{theorem}\cite{yos-3arr}
\label{thm:3arr}
If $\ell=3$, then the cokernel of the restriction map is 
finite dimensional. 
Furthermore, suppose 
that $\exp(\A^H, m^H)=(e_1, e_2)$, then 
$$
\dim_\bC \Coker=b_3(\bC^3\backslash\bigcup_i H_i)-e_1e_2. 
$$
\end{theorem}

\begin{corollary}\cite{yos-3arr} 
\label{cor:3dim}
Suppose $\ell=3$. Then the following conditions are 
equivalent. 
\begin{itemize}
\item $\A$ is free with exponents $(1, e_2, e_3)$. 
\item $\chi(\A, t)=(t-1)(t-e_2)(t-e_3)$ and 
there exists $H\in\A$ such that $\exp(\A^H, m^H)=(e_2, e_3)$. 
\end{itemize}
\end{corollary}

\begin{rem}
Recently a higher dimensional version of 
Corollary \ref{cor:3dim} has been obtained 
by Schulze \cite{shz-free}. 
\end{rem}

The characterization in the case $\ell\geq 4$ is the following. 
\begin{theorem}\cite{yos-char}
\label{thm:char}
Suppose $\ell\geq 4$. Then an arrangement $\A$ is free with 
exponents $(1, e_2, \ldots, e_\ell)$ if and only if 
there exists $H\in\A$ such that 
\begin{itemize}
\item[(a)] 
$(\A^H, m^H)$ is free with exponents 
$(e_2, \ldots, e_\ell)$, and 
\item[(b)] the localization $\A_x=\{H\in\A\mid x\in H\}$ is 
free for any $x\in H\setminus \{0\}$. 
\end{itemize}
\end{theorem}

\subsection{Freeness for multiarrangements}
\label{subsec:char2}

The notion of multiarrangement is a natural 
generalization of simple arrangement. 
For a $2$-dimensional simple arrangement $\A$, 
it is easy to construct explicit basis of $D(\A)$. 
However, for the case of multiarrangements, describing 
an explicit basis for $D(\A, \m)$ is difficult 
even for $\ell=2$. Wakamiko \cite{waka-exp} gave 
an explicit basis for $D(\A, \m)$ with $\ell=2, |\A|=3$. 
Wakefield and Yuzvinsky \cite{wake-yuz} computed the 
exponents for $\ell=2$ and generic $\A$. Both results 
show that the exponents tend to 
$(\lfloor\frac{|\m|}{2}\rfloor, \lceil\frac{|\m|}{2}\rceil)$, 
where $|\m|=\sum_{H\in\A}\m(H)$. 

\begin{rem}
The above mentioned results remind the author 
results of Dolgachev and Kapranov (cf. \S\ref{subsec:sheaf}) 
and Schenck \cite{sch-ele} 
on stability of $\calT_{\bP^n}(-\log(\cup\A))$. 
It seems natural to ask whether for generic $\A$ with 
$\ell=3$, 
$\widetilde{D}(\A, \m)$ is a stable rank $3$ vector bundle 
on $\bP^2$. 
\end{rem}

Recently several results on $D(\A)$ has been 
generalized to multiarrangements. Abe, Terao and 
Wakefield \cite{atw-char} proved that Solomon-Terao's 
formula 
(\ref{eq:stf}) (and (\ref{eq:phi})) 
gives a polynomial $\chi((\A, \m), t)$ for 
any multiarrangement $(\A, \m)$. 
The  polynomial $\chi((\A, \m), t)$ is called the 
characteristic polynomial of a multiarrangement 
$(\A, \m)$, which is a basic tool for proving 
non-freeness for multiarrangements. 

Another important result on free multiarrangements is 
the Addition-Deletion Theorem \cite{atw-emulti}. 
Let $(\A, \m)$ be a multiarrangement. Choose a 
hyperplane $H_0\in\A$ with $\m(H_0)>0$. One can 
associate two multiarrangements to $(\A, \m, H_0)$ 
as follows. 
\begin{itemize}
\item 
The deletion $(\A', m')$: $\A'=\A$ and the multiplicity 
$m':\A'\rightarrow\bZ_{\geq 0}$ is defined by 
$$
m'(H)=\left\{
\begin{array}{ll}
\m(H)&\mbox{ if }H\neq H_0,\\
\m(H)-1&\mbox{ if }H=H_0. 
\end{array}
\right.
$$
\item 
The restriction $(\A'', \m^*)$: $\A''=\{H\cap H_0\mid H\in
\A, H\neq H_0\}$. Let $X\in\A''$. Then $X$ has codimension two. 
Thus the multiarrangement $\A_X=\{H\in\A\mid H\supset X\}$ with 
the multiplicity $\m|_{\A_X}$ is free. We can choose the 
basis $\theta_X, \psi_X, \partial_3, \dots, \partial_\ell$ with 
$\theta_X\notin\alpha_{H_0}\cdot\Der_S$ and 
$\psi_X\in\alpha_{H_0}\cdot\Der_S$. 
Define the multiplicity 
$\m^*:\A''\rightarrow\bZ_{\geq 0}$ 
by $\m^*(X)=\deg\theta_X$. 
\end{itemize}
The following theorem generalizes the classical 
Addition-Deletion theorem \cite{ter-free12} 
to multiarrangements. 

\begin{theorem}
\label{thm:addel}
\cite{atw-emulti} 
With the notations above, any two of the 
following statements imply the third: 
\begin{itemize}
\item[(i)] $(\A, \m)$ is free with exponents 
$(d_1, \dots, d_\ell)$. 
\item[(ii)] $(\A', \m')$ is free with exponents 
$(d_1, \dots, d_\ell-1)$. 
\item[(iii)] $(\A'', \m^*)$ is free with exponents 
$(d_1, \dots, d_{\ell-1})$. 
\end{itemize}
\end{theorem}
Using Theorem \ref{thm:addel}, one can construct 
a lot of free multiarrangements inductively.

\subsection{Free arrangements and intersection of plane curves}
\label{subsec:curv}

In this section we consider a $3$-dimensional 
free arrangement $\A$ with exponents $(1, e_1, e_2)$. 
Choose $H_0\in\A$. Then the deconing ${\bf d}_{H_0}\A$ is 
an affine line arrangement in $\bC^2$. 
Freeness of $\A$ imposes strong conditions on 
the positions of intersections 
$L_2(\A):=\{L\cap L'\in\bC^2\mid L, L'\in{\bf d}_{H_0}\A, 
L\neq L'\}$ and their multiplicities. Let 
$\mu(p):=\sharp\{L\in{\bf d}_{H_0}\A\mid L\ni p\}-1$. 

\begin{theorem}
\label{thm:curv}
Assume that $\A$ is free. 
With notation as above, there exist plane curves 
$C_1, C_2\subset\bC^2$ 
with degrees $e_1$ and $e_2$ respectively 
such that $C_1\cap C_2=L_2(\A)$ and the intersection 
multiplicity is 
$$
\mult_p(C_1, C_2)=\mu(p). 
$$
\end{theorem}

\begin{rem}
If $\A$ is a fiber-type arrangement, we can find easily 
such $C_1$ and $C_2$ as union of lines. 
\end{rem}

\begin{proof}
Choose coordinates $(z_0, z_1, z_2)$ so that 
$H_0=\{z_0=0\}$. We can choose a basis 
$\theta_E, \delta_1, \delta_2\in D(\A)$ 
such that $\delta_1 z_0=\delta_2 z_0=0$ (see 
Lemma \ref{lem:zie}). Let $\alpha=a_1z_1+a_2z_2$ be 
a linear form such that the line $\{\alpha=0\}\subset\bC^2$ 
is not parallel to any line $L\in{\bf d}_{H_0}\A$. 
By definition $\delta_i\alpha\in\bC[z_1, z_2]$ is a polynomial 
of degree $e_i$. Note that 
$$
\delta_i\alpha=0 \mbox{ at }p\in\bC^2
\Longleftrightarrow
\left\{
\begin{array}{l}
\delta_i(p)=0\mbox{ or }\\
\delta_i(p)\mbox{ is parallel to }\{\alpha=0\}. 
\end{array}\right.
$$
If $p\notin\bigcup_{L\in{\bf d}_{H_0}\A}L$, 
$\delta_1(p)$ and $\delta_2(p)$ are linearly independent. 
If $p\in L$ and 
$p\notin L_2(\A)$, %\bigcup_{L'\in{\bf d}_{H_0}\A\setminus\{L\}}L'$, 
$\delta_1(p)$ and $\delta_2(p)$ spans the tangent space 
$T_pL$. In any case, either 
$\delta_1(p)\alpha\neq 0$ or $\delta_2(p)\alpha\neq 0$.  
Hence $\delta_1(p)\alpha=\delta_2(p)\alpha=0$ 
precisely when 
$p\in L_2({\bf d}_{H_0}\A)$. 
Fix $p\in L_2(\A)$ and 
choose the coordinate $(z_1, z_2)$ such that 
$p=(0, 0)$ and $\{z_1=0\}\in{\bf d}_{H_0}\A$. 
Let $Q$ be the product of defining 
equations which contain $p$. Then 
$$
\eta_1=z_1\partial_{z_1}+z_2\partial_{z_2} \mbox{ and }
\eta_2=\frac{Q}{z_1}\partial_{z_2}
$$
form a basis of $D({\bf d}_{H_0}\A_p)$. 
It is easily seen that 
$\mult_p(\eta_1\alpha, \eta_2\alpha)=
\dim \bC[[z_1, z_2]]/(\eta_1\alpha, \eta_2\alpha)=
\mu(p)$. Germs of $\delta_i$ at $p$ is expressed as 
$\delta_i=f_{i1}\eta_1+f_{i2}\eta_2$ with 
$\det=f_{11}f_{22}-f_{12}f_{21}$ is contained in the unit  
$\bC[[z_1, z_2]]^\times$. 
Thus intersection multiplicity is 
$$
\mult_p(C_1, C_2)= 
\dim \bC[[z_1, z_2]]/(\delta_1\alpha, \delta_2\alpha)= 
\mu(p).
$$ 
\end{proof}

%{\bf d}_{H_0}\A

\begin{rem}
Although there exists a free arrangement 
which has non-vanishing homotopy group $\pi_2(M(\A))$, 
\cite{er-h}, 
it is challenging to see the homotopy types of 
free arrangements. 
\end{rem}

\subsection{An example of a non-free arrangement}
\label{subsec:nonfree}

Factorization of the characteristic polynomial 
(Corollary \ref{cor:tfact}) is 
a necessary combinatorial condition for an arrangement to 
be free. However the converse is not true. Indeed, 
there are non-free arrangements 
which have factored characteristic 
polynomials. 

\begin{example}
(Stanley's example \cite{ot}) 
Let $\A=\{H_0, H_1, \dots, H_6\}$ 
be an arrangement of $7$ planes in 
$\bR^3$ defined as Figure \ref{fig:ex} (real lines). 
The characteristic polynomial is $\chi(\A, t)=(t-1)(t-3)^2$. 
However $\A$ is not free. 
We shall give three proofs. 
\begin{figure}[htbp]
\begin{picture}(100,130)(0,0)
\thicklines

\put(30,0){\line(1,1){130}}
\put(40,0){\line(1,2){65}}
\put(50,0){\line(0,1){130}}
\put(60,0){\line(-1,2){65}}
\put(70,0){\line(-1,1){130}}

\put(-48,120){$H_2$}
\put(2,120){$H_3$}
\put(52,120){$H_4$}
\put(107,120){$H_5$}
\put(157,120){$H_6$}

\put(-50,110){\line(1,0){200}}

\put(-50,50){\line(1,0){170}}
\qbezier(120,50)(130,50)(160,20)
\put(160,20){\line(1,-1){10}}
\put(-40,40){$H_1$}

\qbezier(150,110)(160,110)(160,100)

\put(160,100){\line(0,-1){100}}
\put(165,90){$H_0$}

\multiput(-50,20)(5,0){45}{\circle*{2}}
\put(-10,10){$K$}

\end{picture}
\caption{$\A=\{H_0, \dots, H_6\}$}
\label{fig:ex}
\end{figure}
\end{example}

First note that by \ref{cor:tfact}, if $\A$ is free, then the 
exponents should be $(1, 3, 3)$. 

(1) Consider another hyperplane $K$ (dotted line). 
The extended arrangement $\A\cup\{K\}$ is of fiber-type 
and hence free with exponents $(1, 2, 5)$ (also easily proved 
by using Addition-Deletion Theorem \ref{thm:addel}). 
Hence $D(\A\cup\{K\})$ has degree $2$ element $\delta$ which is 
linearly independent from the Euler vector field $\theta_E$. 
By definition, $\delta\in D(\A)$. However this contradicts 
the fact that $D(\A)$ does not have a basis element of 
degree $\leq 2$ other than $\theta_E$. 

(2) Consider the restriction to $H_0$. Then 
$(\A^{H_0}, m^{H_0})$ is free with exponents $(1, 5)$. 
From Corollary \ref{cor:3dim}, $\A$ is not free. 

(3) Consider the deconing ${\bf d}_{H_0}\A$ with respect 
to $H_0$. If $\A$ is free, then by Theorem \ref{thm:curv}, 
the intersections satisfy 
$L({\bf d}_{H_0}\A)=C_1\cap C_2$, where $C_i$ is a 
cubic curve. We may assume that $C_1$ does not have 
$H_1$ as a component. 
Then $H_1\cap C_1$ consists of five points. 
This contradicts Bezout theorem.

\section{Coxeter multiarrangements}
\label{sec:coxeter}

Coxeter multiarrangements are a well-studied 
class of multiarrangements. Using the notion of 
primitive derivation, we can construct a basis 
for several Coxeter multiarrangements. 
Here we give a brief review. 

The importance of the primitive derivation 
was first realized by K. Saito \cite{sai-lin} 
in the context of singularity theory. 
K. Saito's theory of primitive forms reveals 
that the parameter space $B$ of semi-universal 
deformation $X\rightarrow B$ of an isolated 
singularity $0\in X_0$ possesses rich geometric 
structures \cite{sai-per, mats-prim}. 
On the other hand, Grothendieck-Brieskorn-Slodowy's 
theory \cite{bri-sing} shows that for simple singularities, 
the semi-universal family can be described in 
terms of Lie theory. In particular, 
the parameter space $B$ can be canonically 
identified with the Weyl group quotient 
$\frakh/W$ of an $ADE$-type Cartan subalgebra $\frakh$ 
(see also \S\ref{subsec:adj}). 
In \cite{sai-lin}, Saito describes the flat structure for 
any finite reflection group $W\curvearrowright V$ 
in purely invariant 
theoretic way by using the primitive derivation. 
Later Terao \cite{ter-multi, ter-hodge} pulled back the 
theory to $V$ via the natural projection 
$\pi:V\rightarrow V/W$ and proved freeness 
of Coxeter multiarrangements with constant 
multiplicity. 

In this section, we will describe the structure 
of $D(\A, \m)$ for a Coxeter arrangement $\A$ based on 
\cite{ter-multi, ter-hodge, yos-prim, ay-quasi}. 

Let $V$ be an $\ell$-dimensional Euclidean space over $\bR$ 
with inner product $I:V\times V\rightarrow \bR$. 
Fix a coordinate 
$(x_1, \cdots, x_\ell)$ and put 
$S=S(V^*)\otimes_\bR\bC=\bC[x_1, \ldots,x_\ell]$. 
Let $W\subset O(V, I)$ be a finite irreducible reflection group with 
the Coxeter number $h$. 
Let $\A$ be the corresponding Coxeter arrangement, i.e., 
the collection of all reflecting hyperplanes of $W$. 
Fix a defining linear form $\alpha_H\in V^*$ for each 
hyperplane $H\in\A$. 

It is proved by Chevalley \cite{che} that the invariant ring 
$S^W$ is a polynomial ring $S^W=\bC[P_1, \ldots, P_\ell]$ 
with $P_1, \ldots, P_\ell$ are homogeneous generators. 
Suppose that $\deg P_1\leq\cdots\leq\deg P_\ell$. Then it is 
known that 
$\deg P_1=2<\deg P_2\leq\cdots\leq\deg P_{\ell-1}<\deg P_\ell=h$. 
Note that we may choose $P_1(x)=I(x,x)$. 
Then $\frac{\partial}{\partial P_i}$ 
($i=1, \ldots, \ell$) can be 
considered as a rational vector field on $V$ with order one 
poles along $H\in\A$. Indeed by using the fact 
$$
\Delta:=
\det
\left(
\frac{\partial P_i}{\partial x_j}
\right)_{i, j=1,\dots, \ell}
\doteq\prod_{H\in\A}\alpha_H, 
$$
we may define the action of the differential operator 
$\frac{\partial}{\partial P_i}$ to 
$f\in S$ by 
$$
\frac{\partial f}{\partial P_i}=
\frac{1}{\Delta}
\det\left(
\begin{array}{ccccc}
\frac{\partial P_1}{\partial x_1}&
\dots&
\frac{\partial f}{\partial x_1}&
\dots&
\frac{\partial P_\ell}{\partial x_1}\\
\frac{\partial P_1}{\partial x_2}&
\dots&
\frac{\partial f}{\partial x_2}&
\dots&
\frac{\partial P_\ell}{\partial x_2}\\
\vdots&\ddots&\vdots&\ddots&\vdots\\
\frac{\partial P_1}{\partial x_\ell}&
\dots&
\frac{\partial f}{\partial x_\ell}&
\dots&
\frac{\partial P_\ell}{\partial x_\ell}
\end{array}
\right)
$$
Obviously, we have 
$\frac{\partial P_i}{\partial P_i}=1$ and 
$\frac{\partial P_j}{\partial P_i}=0$ for $i\neq j$. 

\begin{defn}
We denote $D=\frac{\partial}{\partial P_\ell}$ and call it 
the {\em primitive derivation}. 
\end{defn}
Since $\deg P_i<\deg P_\ell$ for $i\leq\ell-1$, 
the primitive derivation $D$ is uniquely determined 
up to nonzero constant multiple 
independent of the choice of the generators 
$P_1, \ldots, P_\ell$. 

Next we define the affine connection $\nabla$. 
\begin{defn}
For a given rational vector field 
$\delta=\sum_{i=1}^\ell f_i\frac{\partial}{\partial x_i}$ 
and a rational differential $k$-form 
$\omega=\sum_{i_1, \ldots, i_k}g_{i_1, \ldots, i_k} dx_{i_1, \ldots, i_k}$ 
(where $dx_{i_1, \ldots, i_k}=dx_{i_1}\wedge\dots\wedge dx_{i_k}$), 
define $\nabla\!_\delta\omega$ by 
$$
\nabla\!_\delta\omega=
\sum_{i_1, \ldots, i_k}\delta(g_{i_1, \ldots, i_k}) dx_{i_1, \ldots, i_k}. 
$$
\end{defn}

Let ${\bf m}:\A\longrightarrow\{0, 1\}$ be a map. 
The differentiation $\nabla_D$ by the primitive 
derivation changes the degree by $h$. This action 
connects $D(\A, {\bf m})$ with $D(\A, 2k+{\bf m})$ and 
$\Omega^1(\A, 2k-{\bf m})$. 
\begin{theorem}
\label{thm:ay1}
Fix notation as above, and let $k$ be a positive integer. 
\begin{itemize}
\item[(1)] 
The map 
$$
\begin{array}{cccc}
\Phi_k:&D(\A, \m)(kh)&\longrightarrow&\Omega^1(\A, 2k-\m)\\
&&&\\
&\delta&\longmapsto&\nabla\!_\delta\nabla_D^k dP_1
\end{array}
$$
gives an $S$-isomorphism of graded modules. 
\item[(2)] 
The map 
$$
\begin{array}{cccc}
\Psi_k:&D(\A, \m)(-kh)&\longrightarrow&D(\A, 2k+\m)\\
&&&\\
&\delta&\longmapsto&\nabla\!_\delta\nabla_D^{-k}E
\end{array}
$$
gives an $S$-isomorphism of graded modules. 
\end{itemize}
\end{theorem}

\begin{corollary}
\label{thm:ay2}
For a $\{0,1\}$-valued multiplicity $\m:\A\rightarrow\{0,1\}$ and 
an integer $k>0$, 
the following conditions are equivalent. 
\begin{itemize}
\item $(\A, \m)$ is free with exponents $(e_1, \ldots, e_\ell)$. 
\item $(\A, 2k+\m)$ is free with exponents $(kh+e_1, \ldots, kh+e_\ell)$. 
\item $(\A, 2k-\m)$ is free with exponents $(kh-e_1, \ldots, kh-e_\ell)$. 
\end{itemize}
\end{corollary}
If $\m\equiv 0$, then $(\A, \m)$ is free with 
exponents $(0, \dots, 0)$. Hence $(\A, 2k)$ is 
free with exponents $(kh, kh, \dots, kh)$. 
If $\m\equiv 1$, then $(\A, \m)$ is free with 
exponents $(e_1, \dots, e_\ell)$, where 
$e_i=\deg P_i-1$ (by \cite{sai-log, sai-lin}). 
Hence $(\A, 2k+1)$ is free with exponents 
$(e_1+kh, \dots,  e_\ell+kh)$. 
In particular, Coxeter multiarrangements with 
constant multiplicities are free \cite{ter-multi}. 

The primitive derivation acts on $W$-invariant forms. 
The following will be used in the next section. 

\begin{theorem}
\cite{ter-hodge} 
\label{thm:hfilt}
With notation as above, the set of $W$-invariant 
derivations $D(\A, 2k+1)^W$ is a free $S^W$-module. 
Furthermore, if $k>0$, 
$$
\nabla_{\frac{\partial}{\partial P_i}}D(\A, 2k+1)^W
\subset D(\A, 2k-1)^W, 
\nabla_{\frac{\partial}{\partial P_\ell}}D(\A, 2k+1)^W
=D(\A, 2k-1)^W. 
$$
\end{theorem}

%{\bf m}

\begin{rem}
Recently 
Theorem \ref{thm:ay1} and 
Corollary \ref{thm:ay2} are 
generalized for $\m:\A\longrightarrow\{-1, 0, +1\}$ 
by Abe \cite{abe-gen}. 
\end{rem}

\section{Applications of freeness of $D(\A, \m)$}
\label{sec:appl}

In this section we will describe two applications 
of freeness of $D(\A, \m)$. 

\subsection{Relative de Rham cohomology of adjoint quotient maps}
\label{subsec:adj}

Let $\frakg$ be a simple Lie algebra of type ADE 
over $\bC$. 
The categorical quotient map 
$\chi :\frakg \rightarrow B := \frakg//G $ 
of the adjoint group action 
on $\frakg$ is called the adjoint quotient map. 
The purpose of this section is to investigate the 
$\calO_B$-module  structure of 
the relative de Rham cohomology $H ^2(\Omega_{\chi}^\bullet)$ 
(see below the definition of $\Omega_\chi^\bullet$) 
of $\chi :\frakg \rightarrow B$ through an 
action of vector fields $\Der_B$ on $B$ 
(the Gau\ss -Manin connection).

The study of the relative de Rham cohomology for an 
affine morphism goes back to 
E. Brieskorn \cite{bri-die}  who proved the coherence 
of relative de Rham 
cohomology for any polynomial map $f:\bC ^n \rightarrow \bC$ with 
isolated 
critical point $0\in \bC^n$
(and M. Sebastiani proved $\mathcal{O}_\bC$-freeness of rank $\mu$, 
where $\mu$ is the Milnor number of $f$).
Further, K. Saito proved the freeness for the semi-universal
deformation $F:X\rightarrow B:= \bC ^\mu$ of an isolated hypersurface 
singularity defined by $f$.
More precisely, he gave an isomorphism 
between a certain submodule of vector fields $\Der_B$ on $B$ and 
$H^n(\Omega_{X/B}^\bullet)$. 
The isomorphism is given by the following correspondence,
we first fix a special 
cohomology class $\zeta$ called a primitive form, then 
for given vector field $\delta \in \Der_B$ take a lift up
$\tilde{\delta }\in \Der _X$ of the vector field on the total space $X$,
and differentiate $\zeta$ by $\tilde{\delta}$, we have a new cohomology
class $\mathcal{L}_{\tilde{\delta}}\zeta$, where $\mathcal{L}$
is the Lie derivative. 
On the other hand, the semi-universal deformation of a simple
singularity is constructed by using the 
adjoint quotient map $\chi$ of type ADE 
\cite{bri-sing, slo-simp}. Indeed, 
if we restrict the map $\chi$ 
to a certain affine subspace $X\subset \frakg$, 
we have the semi-universal deformation of a simple singularity. 
In this case H. Yamada \cite{yam-lie} showed that the restriction 
of the Kostant-Kirillov form $\zeta$ to $X$ becomes 
the primitive form 
which generates the relative de Rham cohomology 
$H^2(\Omega_{X/B}^\bullet)$ 
by differentiation by means of vector fields $\delta \in \Der_B$.

Let us introduce some notation. 
Let $\frakg$ a simple Lie algebra over $\bC$ 
(later we will restrict $\frakg$ to ADE-type). 
Let $\frakg =\frakh \oplus \bigoplus_{\alpha \in \Phi} 
\frakg_{\alpha}, 
(\Phi \subset \frakh^{\ast})$ a Cartan decomposition 
with respect to a Cartan algebra $\frakh$ with 
$\ell =\dim \frakh$, $G$ the adjoint group of 
$\frakg$ and $T$ the maximal torus of $G$ 
with Lie algebra $\frakh$. 
We denote by $W$ the Weyl group $N_G (T) /T$ .
The classical Chevalley's restriction theorem states that 
the restriction $\rho :\bC [\frakg]\rightarrow \bC [\frakh]$ 
of polynomial functions induces an isomorphism 
\begin{equation}
\label{eq:rest1}
\bC [\frakg]^G \stackrel{\cong}{\longrightarrow} \bC [\frakh]^W
\end{equation}
of algebras of invariants. We also denote 
$\bC [\frakh]=S$ and $S^W=\bC[P_1, \dots, P_\ell]$ as 
in \S\ref{sec:coxeter}. The categorical quotient of 
the adjoint action is 
$B= \frakg //G \cong \frakh//W\cong \Spec S^W$. 
We call the quotient map $\chi:\frakg\rightarrow B$ 
the adjoint quotient map as mentioned above. 
The construction is summarized in the following 
diagram. 
\begin{equation}
\label{eq:diag}
\begin{array}{ccccccc}
&&\frakg&&&&\bC[\frakg]\\
&&\downarrow&&&&\uparrow\\
V=\frakh&\stackrel{\pi}{\longrightarrow}&\frakh/W=B=\frakg//G&, 
&S&\hookleftarrow&S^W=\bC[B]=\bC[\frakg]^G. 
\end{array}
\end{equation}
\begin{defn}
Define the relative de Rham complex $\Omega_\chi^\bullet$ 
for the adjoint quotient map $\chi:\frakg\rightarrow B$ by 
$$
\Omega_\chi^\bullet=
\frac{\Omega_\frakg^\bullet}{\chi^*\Omega_B^1\wedge
\Omega_\frakg^{\bullet-1}}
=
\frac{\Omega_\frakg^\bullet}{\sum_{i=1}^\ell dP_i\wedge
\Omega_\frakg^{\bullet-1}}. 
$$
\end{defn}
By the formula $d(P\cdot\omega)=dP\wedge\omega+P\cdot d\omega$, 
the differential 
$d_\chi:\Omega_\chi^\bullet\rightarrow\Omega_\chi^{\bullet+1}$ 
is a $S^W$-module homomorphism. Hence the cohomology 
group 
$H^k(\Omega_\chi^\bullet)$ possesses $S^W$-module structure. 

%Let $\A$ be the set of reflecting hyperplanes in $\frakh$. 
Let $D\subset B$ be the set of critical points of the 
quotient map $\pi:\frakh\rightarrow B$. It is 
proved in \cite{sai-log} that $\pi$ induces an 
isomorphism 
$$
D(\A)^W\stackrel{\cong}{\longrightarrow} \Der_B(-\log D). 
$$
Thus for $\delta\in D(\A)^W$, we may differentiate 
the Kostant-Kirillov form $\zeta$ by 
$\delta$ and obtain a relative $2$-form 
$\nabla_\delta\zeta$ (which has poles 
along $D$ in general). 

After Yamada's result, it was naturally 
conjectured that $H^2(\Omega_\chi^\bullet)$ is 
a free $S^W$-module of rank $\ell$. 

\begin{theorem}
\label{thm:adjqt}
Let $\frakg$ a simple Lie algebra of type ADE, 
with a Cartan subalgebra $\frakh$ and the Weyl group 
$W$. Let $\A$ be the corresponding Weyl arrangement 
on $\frakh$. The map $D(\A)^W\ni\delta
\longmapsto \nabla_\delta\zeta$ induces 
a natural isomorphism 
$$
H^2(\Omega_{{\chi}}^\bullet)
\cong
D(\A, 5)^W, 
$$
of $S^W$-modules. 
\end{theorem}
The rest of this section is devoted to a proof of 
this theorem. 

We first recall a result due to J. Vey \cite{vey}, 
which is an analogue of Weyl's unitary trick. 

\begin{theorem}
\label{thm:vey}
Let $G$ be a connected reductive algebraic group over $\bC$ 
with a linear action on a finite dimensional $\bC$-vector space $E$. 
Let  $\Omega_E^\bullet$
be the de Rham complex of holomorphic differential forms on $E$ 
and $\mathcal{I}^\bullet$ the ideal of $\Omega^\bullet$ generated by 
differentials $df_1, df_2, \cdots ,df_r$, where 
$f_1, f_2, \cdots ,f_r$ 
are $G$-invariant homogeneous polynomials on $E$. Then the 
morphism 
$$
(\Omega ^\bullet )^G /(\mathcal{I} ^\bullet )^G \rightarrow 
\Omega^\bullet  /\mathcal{I}^\bullet 
$$
is a 
quasi-isomorphism. 
\end{theorem}
By this, we can compute cohomology of $\Omega_\chi^\bullet$ 
by using the complex 
$\Omega_\frakg^{G, \bullet}/(\sum_i dP_i\wedge
\Omega_\frakg^{\bullet-1})^G$ of 
$G$-invariant relative forms. 
Next we shall describe 
$\Omega_\frakg^{G, \bullet}$. 
%Using Broer's result. 

Broer \cite{bro} considered a generalization of Chevalley's 
restriction theorem (\ref{eq:rest1}) in the following setting. 
Let $M$ be a finite dimensional $G$-module 
and $\Mor (\frakg , M)$ (resp. $\Mor_G (\frakg , M)$) the space 
of polynomial (resp. $G$-equivariant polynomial) morphisms 
of $\frakg$ into $M$. It is isomorphic to $\bC [\frakg]\otimes M$ 
(resp. $(\bC [\frakg]\otimes M)^G$). For any $G$-module $M$ 
the restriction map $\rho$ induces a homomorphism 
\begin{displaymath}
\rho_M :\Mor _G(\frakg, M) \longrightarrow \Mor_W(\frakh ,M^T). 
\end{displaymath}
Since the union of all Cartan subalgebras is Zariski 
dense in $\frakg$, 
$\rho_M$ is injective for all $M$. 
If $M=\bC$ is a trivial $G$-module, $\rho_M$ is bijective because of 
Chevalley's theorem. 
However it is not necessarily bijective in general. 
Broer \cite{bro} proved that 

\begin{theorem}\label{thm:bro}
Let $M$ be a $G$-module. Restriction induces an isomorphism 
\begin{displaymath}
\rho_M :\Mor _G\left(\frakg ,M\right) \stackrel{\cong}{\longrightarrow} 
\Mor_W\left(\frakh ,M^T\right) 
\end{displaymath}
if and only if the weights $2\alpha$ 
($\alpha \in \Phi$ is a root 
of $\frakg$) do not occur as $T$-weights in $M$. 
(We shall call $M$ {\em small} if it satisfies this assumption.)
\end{theorem}
We need this theorem to describe the set of $G$-invariant 
differential forms $\Omega_\frakg^{\bullet,G}$ on $\frakg$ below. 
By definition the set of all differential $p$-forms on $\frakg$ is 
$\Omega_\frakg^p =\bC [\frakg]\otimes \stackrel{p}{\wedge}\frakg^\ast$. 
Thus we apply Theorem \ref{thm:bro} for 
$M= \stackrel{p}{\wedge}\frakg^\ast$. 
If $p=1$, since $\frakg\cong \frakg^\ast $ by Killing form, the 
$T$-weights of $\frakg^\ast$ are nothing but the roots of $\frakg$, 
so $\frakg^\ast $ is small. Thus we have an isomorphism 
\begin{displaymath}
\rho_1 :\Omega_\frakg^{1,G} \cong \left(\bC [\frakh]\otimes \frakh^\ast 
\right)^W
\cong \Omega_\frakh^{1,W}.
\end{displaymath}
It follows from a result of Solomon \cite{sol-inv} that 
\begin{equation}
\label{eq:sol}
\pi^*:\Omega_B^p\stackrel{\cong}{\longrightarrow}
\Omega_\frakh^{p,W}. 
\end{equation}
Thus we conclude that 
$G$-invariant $1$-forms $\Omega_\frakg^{1,G}$ on 
$\frakg$ are nothing but pull back $\chi^\ast \Omega_B^1$ of 
$1$-forms 
on $B$. In particular, $\Omega_{{\chi}}^{1,G} =0$, 
we have 
\begin{equation}\label{eq:ker}
H^2(\Omega_{{\chi}}^{\bullet ,G})\cong \Ker 
(d_{{\chi}}: \Omega_{{\chi}}^{2,G}\rightarrow 
\Omega_{{\chi}}^{3,G}). 
\end{equation}
From the classification of simple root systems, 
it is easily seen that $M=\stackrel{2}{\wedge}\frakg^\ast$ is small 
if and only if $\frakg$ is of type ADE, since the set of weights 
of $\stackrel{2}{\wedge}\frakg^\ast$ is 
\begin{displaymath}
\{ 0\} \cup \Phi \cup 
\{ \alpha +\beta | \alpha , \beta \in \Phi ,\alpha \neq \beta \}.
\end{displaymath}

Furthermore, 
$(\stackrel{2}{\wedge}\frakg^\ast )^T 
\cong \stackrel{2}{\wedge}\frakh^\ast\otimes 
(\stackrel{2}{\wedge}\frakh^\perp)^T$ is a direct sum decomposition of 
$W$-submodules. From (\ref{eq:sol}), we obtain, 
\begin{prop}
\label{prop:decomp}
Let $\frakg$ be a simple Lie algebra of type ADE. Then 
\begin{equation}
\label{eq:directsum}
\Omega_\frakg^{2,G}\cong 
\chi^\ast \Omega_B^2 \oplus \Omega_{{\chi}}^{2,G}. 
\end{equation}
\begin{equation}
\label{eq:rel2form}
\rho_2:
\Omega_{{\chi}}^{2,G} \stackrel{\cong }{\longrightarrow}
\left(\bC[\frakh]\otimes 
(\stackrel{2}{\wedge}{\frakh^\perp})^T\right)^W. 
\end{equation}
\end{prop}

By Proposition \ref{prop:decomp}, 
we can identify $G$-invariant relative $2$-forms 
$\Omega_{{\chi}}^{2,G}$ with the submodule 
$(\bC[\frakh]\otimes (\stackrel{2}{\wedge}{\frakh^\perp})^T)^W$ 
of $\Omega_\frakg^{2,G}$. 
Now we define two submodules 
$\calH_\chi\subset\calH_\chi'\subset \Omega_\frakg^{2,G}$ 
which are related to $H^2(\Omega_\chi^\bullet)$. 

\begin{defn}

\begin{equation}
\label{eq:hchi}
\calH_\chi:=
\{ \omega \in \Omega_{{\chi}}^{2,G} 
\mid 
d_\frakg\omega\in \sum_i dP_i\wedge\Omega_\frakg^{2,G} \}, 
\end{equation}
\begin{equation}
\label{eq:hchi'}
\calH_\chi':=
\{ \omega \in \Omega_{{\chi}}^{2,G}\  
\mid 
d_\frakg\omega 
\wedge dP_1 \wedge \cdots \wedge dP_\ell =0 \mbox{\ in\ } 
\Omega_\frakg^{3+\ell}\}, 
\end{equation}
where 
$\bC [\frakg]^G =\bC [P_1, \cdots ,P_\ell ]$. 
\end{defn}
By (\ref{eq:ker}) and definition above, 
\begin{equation}\label{eq:realize}
\calH_\chi\cong H^2(\Omega_{{\chi}}^{\bullet ,G}) 
\end{equation}
and obviously 
$\calH_\chi\subset\calH_\chi'$. 
Later it will be proved that 
$\calH_\chi\subsetneqq \calH_\chi'$. 

Let $e_\alpha \in \frakg_\alpha \ (\alpha \in \Phi )$ be non-zero root 
vectors such that $I( [e_\alpha ,e_{-\alpha}],h)=\alpha (h)$ (for all 
$\alpha \in \Phi ,\ h\in \frakh$), where 
$I(\bullet ,\bullet )$ is the 
Killing form, and $e_\alpha^\ast \in \frakg_\alpha^\ast$ 
be the dual basis. 
Then each element of 
$(\bC[\frakh]\otimes (\stackrel{2}{\wedge}{\frakh^\perp})^T)^W$ 
can be expressed in the form 
\begin{displaymath}
\omega =\sum_{\alpha \in \Phi ^+} f_{\alpha} 
\cdot e^\ast_{\alpha}\wedge 
e^\ast_{-\alpha}
\in \left(\bC[\frakh^{\ast}]\otimes 
(\stackrel{2}{\wedge}{\frakh^\perp})^T\right)^W.
\end{displaymath}
Since $\omega$ is $W$-invariant, if we apply the simple reflection 
$s_\alpha \in W$ with respect to a root $\alpha \in \Phi^+$, 
to $\omega$ 
we have $s_\alpha f_\alpha =-f_\alpha$. Hence $f_\alpha$ is 
divisible by $\alpha$.

Next let us recall the definition of the Kostant-Kirillov form. 
The Kostant-Kirillov form $\zeta$ is a symplectic form 
on the (co)adjoint orbit $G\cdot x\subset\frakg$ of 
$x\in\frakg$. Let $Y, Z\in\frakg$. Then 
$[Y, x]=\left.\frac{d}{dt}\right|_{t=0}\ad(e^{tX})x
\in T_x (G\cdot x)$. For two tangent vectors 
$[Y, x], [Z, x]\in T_x (G\cdot x)$, the $2$-form 
$\zeta$ is given by the formula 
$$
\zeta([Y, x], [Z, x])=I(x, [Y, Z]), 
$$
where $[Y, Z]$ is the bracket in $\frakg$. 

\begin{prop}
\label{prop:kk}
By restricting the Kostant-Kirillov form $\zeta$ to $\frakh$, 
we have 
the following expression 
\begin{equation}
\label{eq:restkk}
\rho_2 (\zeta )=-\sum_{\alpha \in \Phi ^+}
\frac{e^\ast_{\alpha}\wedge e^\ast_{-\alpha}}{\alpha}.
\end{equation}
\end{prop}
\begin{proof}
Let $h\in\frakh\setminus\bigcup H_\alpha$. 
We compute 
$\zeta([e_\alpha, h], [e_\beta, h])$ in two ways. 
First, using the property 
$[h, e_\alpha]=\alpha(h)e_\alpha$, 
we have 
$\zeta([e_\alpha, h], [e_\beta, h])=\alpha(h)\beta(h)
\zeta(e_\alpha, e_\beta)$. 
On the other hand, using the definition of $\zeta$, 
we have 
$\zeta([e_\alpha, h], [e_\beta, h])=I(h, [e_\alpha, e_\beta])$. 
Note that it is non-zero only if $\beta=-\alpha$, and in this case, 
we have $I(h, [e_\alpha, e_{-\alpha}])=\alpha(h)$. 
Hence we have $\zeta(e_\alpha, e_{-\alpha})=\frac{-1}{\alpha(h)}$, 
which implies (\ref{eq:restkk}). 
\end{proof}

The generic fiber of $\chi:\frakg\rightarrow B$ is 
isomorphic to $G/T$, which is homotopy equivalent to the 
flag manifold of $G$. We recall the 
Borel-Hirzebruch description of the de Rham 
cohomology of $G/T$ in degree $2$ \cite{bh}. A $G$-invariant 
differential form on $G/T$ can be seen as a $G$-invariant 
section of the vector bundle 
$\stackrel{\bullet}{\wedge}T^\ast (G/T)$. Hence the evaluation 
at the base point $[T]\in G/T$ induces an isomorphism 
\begin{equation}
\lambda: \Omega_{G/T}^{\bullet,G}
\stackrel{\cong}{\longrightarrow}
\left(\stackrel{\bullet}{\wedge}T_{[T]}^\ast (G/T)\right)^T
\cong 
\left(\stackrel{\bullet}{\wedge}\frakh^\perp \right)^T. 
\end{equation}
For degree $2$, the above map induces the isomorphism 
$\Omega_{G/T}^{\bullet,G}\cong
(\wedge^2\frakh^\perp )^T=
\bigoplus_{\alpha\in\Phi^+}(\frakg_\alpha\oplus\frakg_{-\alpha})$. 
Using the map $\lambda$, we can show that 
$H^2(G/T, \bC)\cong\frakh$. 
\begin{theorem}
\label{thm:bh}
\cite{bh} 
\begin{displaymath}
\begin{array}{cccc}
\omega :&\frakh&\longrightarrow &\Omega_{(G/T)}^{2,G}\\
&&&\\
&h&\longmapsto&
\sum\limits_{\alpha\in\Phi}
\alpha(h)\cdot\lambda^{-1}(e_\alpha^\ast\wedge e_{-\alpha}^\ast)
\end{array}
\end{displaymath}
induces an isomorphism of $\bC$-vector spaces 
$\frakh\stackrel{\cong}{\longrightarrow}
H^2(G/T,\bC)$. 
In particular, $\omega(h)$ is closed, 
\begin{equation}
\label{eq:cl}
\sum_{\alpha\in\Phi}
\alpha(h)\cdot d \lambda^{-1}(e_\alpha^\ast\wedge e_{-\alpha}^\ast)
=0. 
\end{equation}
\end{theorem}

Next we consider the relative de Rham cohomology 
for the projection 
$$
\pro:\frakh\times G/T\rightarrow \frakh. 
$$
We may consider $G$ acts on each fiber of $\pro$ 
from the left. 
Since 
$\Omega_{{\pro}}^2\cong
\bC[\frakh]\otimes\Omega_{G/T}^2$, 
the set of $G$-invariant 
relative $2$-forms $\Omega_{{\pro}}^2$ 
is described by the following isomorphism 
\begin{equation}
\label{eq:lambda}
1\otimes 
\lambda:
\Omega_{{\pro}}^2
\stackrel{\cong}{\longrightarrow}
\bC[\frakh]\otimes (\stackrel{2}{\wedge}\frakh^\perp )^T. 
\end{equation}

By definition and Theorem \ref{thm:bh}, 
$$
H^2(\Omega_{\pro}^\bullet)
\cong
\bC[\frakh]\otimes H^2(\Omega_{(G/T)}^\bullet)
\cong\bC[\frakh]\otimes\frakh 
\cong\Der_\frakh.
$$
The isomorphism is given by

\begin{theorem}
\label{thm:trivial}
\begin{equation}
\label{eq:omega}
\begin{array}{cccc}
1\otimes\omega : & \Der_\frakh & \longrightarrow & 
\bC[\frakh]\otimes
\Omega_{(G/T)}^{2,G}\cong 
\Omega_{\pro}^{2,G}\\
&&&\\
& \delta &\longmapsto & 
(1\otimes\omega)\delta=
\sum\limits_{\alpha\in\Phi}(\delta\alpha)\cdot
\lambda^{-1}(e_\alpha^\ast\wedge e_{-\alpha}^\ast)
\end{array}
\end{equation}
induces a $\bC[\frakh]$-module isomorphism 
$\Der_\frakh\cong 
H^2(\Omega_{\pro}^\bullet)$, 
where 
$\delta\alpha$ is the differentiation of a function 
$\alpha$ by a vector 
field $\delta$. 
\end{theorem}

We define a submodule 
${\calH}_{\pro}\subset \Omega_{(G/T)\times\frakh}^{2,G}$ as 
$$
{\calH}_{\pro} :=
\left\{\left. 
\sum\limits_{\alpha\in\Phi}(\delta\alpha)\cdot 
\lambda^{-1}(e_\alpha^\ast\wedge e_{-\alpha}^\ast)\in 
\Omega_{(G/T)\times\frakh}^{2,G} 
\right|\ \delta \in \Der_\frakh \right\}. 
$$

From the decomposition $\Omega_{(G/T)\times\frakh}^{2,G}=
\bigoplus\limits_{i+j=2}\Omega_{G/T}^{i,G}\wedge\Omega_\frakh^j$ 
and 
$\Omega_{{\pro}}^{2,G}\cong 
\bC[\frakh]\otimes\Omega_{G/T}^{2,G}$, 
we may consider 
\begin{equation}
\label{eq:subset}
\Omega_{(G/T)\times\frakh}^{2,G}
\supset\Omega_{{\pro}}^{2,G}
\supset\calH_{\pro}\ \left(\cong 
H^2(\Omega_{\pro}^\bullet)\right). 
\end{equation}
Let $x_1,\cdots ,x_\ell$ be a coordinate system of $\frakh$, 
then 
$\calH_{\pro}$ 
has another expression as in (\ref{eq:hchi'}): 
\begin{equation}
\calH_{\pro} =
\left\{\left. \omega\in 
\Omega_{{\pro}}^{2,G}\right|\ 
d\omega \wedge dx_1\wedge\cdots\wedge dx_\ell =0 
\mbox{ in } \Omega_{(G/T\times \frakh)}^{3+\ell,G}\right\}. 
\end{equation}

To study the relative de Rham cohomology of 
$\chi:\frakg\rightarrow B$, 
we consider the following diagram 
\begin{equation}
\label{diag:key}
\begin{array}{cccc}
(G/T)\times\frakh &\stackrel{\widetilde{\pi}}{\longrightarrow} 
&\frakg & \\
\pro\downarrow \ \ \ & &\ \ \ \downarrow \chi& \\
\frakh &\stackrel{\pi}{\longrightarrow} &B &
\end{array}
\ \ 
\left(
\begin{array}{cccc}
(g[T],h) &\longmapsto &\ad(g)h & \\
\downarrow & &\downarrow & \\
h &\longmapsto &\overline{h} & 
\end{array}
\right).
\end{equation}
More precisely, from diagram (\ref{diag:key}), there is a natural 
homomorphism 
\begin{displaymath}
\widetilde{\pi}^\ast: 
H^2\left(\Omega_{{\chi}}^\bullet\right)\hookrightarrow 
H^2\left(\Omega_{{\pro}}^\bullet\right), 
\end{displaymath}
which is injective because we have realized 
these cohomology groups 
as subspaces of absolute differential forms 
(see (\ref{eq:realize}) and (\ref{eq:subset})). 
We consider the image of 
$H^2(\Omega_{{\chi}}^\bullet)$ in 
$H^2(\Omega_{{\pro}}^\bullet)\cong \Der_\frakh$. 
Note that if we define a $W$ action on $G/T\times\frakh$ by 
\begin{displaymath}
w\cdot(g[T],h)=(gn_w^{-1},\ad (n_w)h), 
\end{displaymath}
where $n_w\in N_G(T)$ is a representative of $w\in W= N_G(T)/T$, 
then obviously $\widetilde{\pi}$ is a $W$-invariant map and the 
pull back of differential form 
(resp. relative cohomology class) on $\frakg$ 
by $\widetilde{\pi}$ becomes a  
$W$-invariant differential form (resp. $W$-invariant cohomology class). 
\begin{equation}
\widetilde{\pi}^\ast\calH_{\chi} \subset \calH_{\pro}^W. 
\end{equation}

Now recall two expressions of relative $2$-forms 
(\ref{eq:rel2form}) in 
Proposition \ref{prop:decomp} and 
(\ref{eq:lambda}), we have a diagram: 

\begin{equation}
\label{eq:confuse}
\begin{array}{ccccc}
\Omega_{{\chi}}^{2,G} 
&\stackrel{\widetilde{\pi}^\ast}{\hookrightarrow} 
&\Omega_{{\pro}}^{2,G\times W}
&\subset  &\Omega_{{\pro}}^{2,G}\\
\rho_2 \downarrow  \wr & &1\otimes \lambda \downarrow \wr & 
&\downarrow \wr 
\\
\left(\bC [\frakh^\ast]\otimes (\stackrel{2}{\wedge} \frakh^\perp)^T 
\right)^W
&\stackrel{(1\otimes\lambda) \circ 
\widetilde{\pi}^\ast \circ 
\rho_2^{-1}}{\longrightarrow}&
\left(\bC[\frakh ^\ast]\otimes (\stackrel{2}{\wedge} \frakh ^\perp)^T 
\right)^W &\subset &
\bC[\frakh ^\ast]\otimes (\stackrel{2}{\wedge} \frakh ^\perp)^T.
\end{array}
\end{equation}
We compute the map 
$(1\otimes\lambda) \circ \widetilde{\pi}^\ast \circ \rho_2^{-1}$. 

\begin{lem}
\label{lem:bibun}
The map 
$(1\otimes\lambda) \circ \widetilde{\pi}^\ast \circ \rho_2^{-1}: 
\left(\bC[\frakh^\ast]\otimes (\stackrel{2}{\wedge} \frakh^\perp)^T 
\right)^W
\rightarrow
\left(\bC[\frakh ^\ast]\otimes 
(\stackrel{2}{\wedge} \frakh ^\perp)^T 
\right)^W$
can be expressed as 
\begin{displaymath}
\sum\limits_{\alpha \in \Phi^+} \alpha f_\alpha \cdot
e^\ast_\alpha \wedge e^\ast_{-\alpha}
\longmapsto
-\sum\limits_{\alpha \in \Phi^+} \alpha ^3 f_\alpha \cdot
e^\ast_\alpha \wedge e^\ast_{-\alpha}.
\end{displaymath}
\end{lem}
\begin{proof}
The derivation of $\widetilde{\pi}$ is given by 

\begin{displaymath}
\begin{array}{cccc}
(d\widetilde{\pi} )_{([T],h)}:&T_{([T],h)}(G/T\times\frakh )
&\longrightarrow&T_h\frakg\\
 &||\wr & &||\wr\\
 &\frakg/\frakh \oplus \frakh&&\frakg\\
&&&\\
& (\bar{X_1},X_2)&\longmapsto &[X_1,h]+X_2. 
\end{array}
\end{displaymath}
Indeed 
\begin{displaymath}
(d\widetilde{\pi} )_{([T],h)}(\bar{X_1},X_2)=
\left.\frac{d}{dt}\right| _{t=0}
\ad (\exp(tX_1))(h+tX_2)=[X_1,h]+X_2, 
\end{displaymath}
In particular 
$\frakg/\frakh \oplus \frakh \ni (\bar{e}_\alpha ,0) 
\longmapsto
-\alpha (h) e_\alpha \in \frakg _\alpha$. Hence we have 

\begin{displaymath}
\widetilde{\pi}^\ast (e_\alpha^\ast) =-\alpha (h) e_\alpha^\ast \ ,\ 
\widetilde{\pi}^\ast (e_{-\alpha}^\ast) =\alpha (h) e_{-\alpha}^\ast.
\end{displaymath}

\end{proof}

\begin{example}
\label{ex:kk2}
From Proposition \ref{prop:kk} and the preceding lemma, 
the pull back of the 
Kostant-Kirillov form $\zeta$ is 
\begin{displaymath}
\widetilde{\pi}^\ast(\zeta)=
\sum\limits_{\alpha \in \Phi^+} \alpha \cdot
\lambda^{-1}(e^\ast_\alpha \wedge e^\ast_{-\alpha})
\end{displaymath}
Using the Euler vector field 
$\theta_E:=\sum_{i=1}^\ell x_i \frac{\partial}{\partial x_i}$, 
it is expressed as 
$\widetilde{\pi}^\ast(\zeta)=(1\otimes \omega)(\theta_E)$. 
\end{example}

As a corollary of Lemma \ref{lem:bibun}, 
we can characterize the image of the map 
$\widetilde{\pi}^\ast:\Omega_{{\chi}}^{2,G}\rightarrow 
\Omega_{{\pro}}^{2,G\times W}$. 

\begin{corollary}
\label{cor:image}
\begin{displaymath}
(1\otimes\lambda)\circ
\widetilde{\pi}^\ast(\Omega_{{\chi}}^{2,G})
=
\left\{\left.\sum\limits_{\alpha \in \Phi^+} F_\alpha \cdot
e^\ast_\alpha \wedge e^\ast_{-\alpha} \ \right|\ F_\alpha
\mbox{ can be divisible by }\alpha ^3 \right\}^W. 
\end{displaymath}
\end{corollary}
\begin{proof}
If $\sum_{\alpha \in \Phi^+} F_\alpha \cdot e^\ast_\alpha \wedge 
e^\ast_{-\alpha}$
is contained in the left hand side above, $F_\alpha$ have to be 
divisible by $\alpha^3$ from the preceding lemma. Conversely, 
if $F_\alpha$ is divisible by $\alpha^3$ for all $\alpha\in\Phi^+$, 
it is the image of 
\begin{displaymath}
\rho_2^{-1}\left(\sum_{\alpha \in \Phi^+} 
\frac{F_\alpha}{\alpha^2} \cdot 
e^\ast_\alpha \wedge e^\ast_{-\alpha}\right)
\in \Omega_{{\chi}}^{2,G}.
\end{displaymath}
\end{proof}

%%\left(\bbC[\frakh ^\ast]\otimes (\stackrel{2}{\wedge} \frakh ^\perp)^T 
%\right)^W
%%&\stackrel{(1\otimes\lambda) \circ \widetilde{\pi}^\ast \circ 
%\rho_2^{-1}}{\longrightarrow}&
%%\left(\bbC[\frakh ^\ast]\otimes (\stackrel{2}{\wedge} \frakh ^\perp)^T 
%\right)^W

Here it is possible to characterize the image of 
$\widetilde{\pi}^\ast$. 
We have a diagram deduced from (\ref{eq:confuse}), 
\begin{equation}
\label{eq:simple}
\begin{array}{ccrcc}
\calH_{\chi}\subset\calH_{\chi}' 
&\stackrel{\widetilde{\pi}^\ast}{\hookrightarrow} 
&\calH_{\pro}^W&\subset  &\calH_{\pro}\\
& &(1\otimes\omega)^{-1} \downarrow \wr & &\downarrow \wr \\
&&\Der_\frakh^W &\subset &\Der_\frakh.
\end{array}
\end{equation}

Combining (\ref{eq:omega}) and 
Corollary \ref{cor:image}, we have 

\begin{theorem}
\label{thm:h'}
$(1\otimes\omega)^{-1}\circ\widetilde{\pi}^\ast$ 
induces an isomorphism 
\begin{displaymath}
\calH_{\chi}'\cong D(\A, 3)^W. 
\end{displaymath}
\end{theorem}
\begin{proof}
For any $\delta\in\Der_\frakh^W$, since 
\begin{eqnarray*}
d\left((1\otimes\omega)\delta\right)\wedge 
dP_1\wedge\cdots\wedge dP_\ell  
&=&
Q\cdot d \left((1\otimes\omega)\delta\right)\wedge 
dx_1\wedge\cdots\wedge 
dx_\ell\\
&=& 0,
\end{eqnarray*}
$(1\otimes\omega)\delta \in 
\widetilde{\pi}^\ast \calH_\chi'$ if and only if 
$\delta\alpha$ is divisible by $\alpha^3$ for all 
$\alpha \in \Phi^+$.
\end{proof}

Now we are in a position to prove our main result. 

\begin{theorem}
\label{thm:h}
$(1\otimes\omega)^{-1}\circ\widetilde{\pi}^\ast$ 
induces an isomorphism 
\begin{displaymath}
\calH_{\chi}\cong D(\A, 5)^W. 
\end{displaymath}
Hence $H^2(\Omega_{{\chi}}^\bullet)$ is a free 
$\bC[B]$-module of rank $\ell$. 
\end{theorem}
\begin{proof}
Suppose $\eta\in\calH_{\chi}$ and put 
$\widetilde{\pi}^\ast \eta 
=\sum_{\alpha\in\Phi}
(\delta\alpha)\cdot
\lambda^{-1}(e_\alpha^\ast\wedge e_{-\alpha}^\ast)$ 
for $\delta\in D(\A, 3)^W$, 
then by definition there exist 
$\eta_1, \cdots,\eta_\ell\in \Omega_{\frakg}^{2,G}$ such that 

\begin{displaymath}
d\eta =\sum_{i=1}^{\ell}dP_i\wedge\eta_i.
\end{displaymath}
Applying the operator $d$ and multiplying by 
$dP_1\wedge\cdots \wedge\widehat{dP_i}\wedge\cdots\wedge dP_\ell$, 
\begin{displaymath}
dP_1\wedge\cdots\wedge dP_i\wedge\cdots\wedge 
dP_\ell \wedge d\eta_i =0.
\end{displaymath}
Hence $\eta_i \in \calH_{\chi}'$ for all $i=1,\cdots,\ell$ and 
$\widetilde{\pi}^\ast\eta_i 
=\sum_{\alpha\in\Phi}
(\delta_i\alpha)\cdot\lambda^{-1}
(e_\alpha^\ast\wedge e_{-\alpha}^\ast)$
for some $\delta_i \in D(\A, 3)^W$. 
Using (\ref{eq:cl}), we have 
\begin{eqnarray*}
d\widetilde{\pi}^\ast \eta &=& 
\sum_{\alpha\in\Phi}d(\delta\alpha)\wedge 
\lambda^{-1}(e_\alpha^\ast\wedge e_{-\alpha}^\ast) \\
&=& \sum_{\alpha\in\Phi}\sum_i
\frac{\partial}{\partial P_i}(\delta\alpha)dP_i \wedge 
\lambda^{-1}(e_\alpha^\ast\wedge e_{-\alpha}^\ast) \\
&=&
\sum_{\alpha\in\Phi}\sum_i
\left((\nabla_\frac{\partial}{\partial P_i} \delta)\alpha\right)
dP_i \wedge 
\lambda^{-1}(e_\alpha^\ast\wedge e_{-\alpha}^\ast), 
\end{eqnarray*}
and $\delta_i =\nabla_\frac{\partial}{\partial P_i} \delta$. 
By Theorem \ref{thm:hfilt}, $\eta \in \calH_{\chi}$ if and only if 
$\delta \in D(\A, 5)^W$. 
\end{proof}

\subsection{Freeness of $A_n$-Catalan arrangements}
\label{subsec:cat}

As another application, we prove that 
Catalan arrangements of type $A$ are free. 

Let $(x_1, \dots, x_n, z)$ be a coordinate of 
$\bC^{n+1}$. 
(The cone of) Catalan arrangement $\Cat_n$ 
is defined by 
$$
z\times\prod_{1\leq i<j\leq n}
\left(
(x_i-x_j)(x_i-x_j-z)(x_i-x_j+z)
\right)=0. 
$$
The terminology ``Catalan arrangement'' comes from the 
fact that the number of chambers divided by $2n!$ is 
equal to $n$-th Catalan number. See \cite{ath-survey} 
for more combinatorial aspects of $\Cat_n$. (Note 
that the definition of $\Cat_n$ in this article is 
the coning of that of \cite{ath-survey}.) 

\begin{theorem}
\label{thm:cat}
The Catalan arrangement $\Cat_n$ is free with 
exponents $(0, 1, n+1, n+2, \dots, 2n-1)$. 
\end{theorem}

\begin{rem}
This result was first proved by 
Edelman and Reiner \cite{ede-rei} using 
Addition-Deletion Theorem \ref{thm:addel}. 
It can be also proved by using 
the freeness of $D(\A, 3)$ and Theorem 
\ref{thm:char}. We give another proof 
which is also based on the freeness of 
$D(\A, 3)$. However, instead of using 
Theorem \ref{thm:char}, 
we will directly show the existence of basis of 
$D(\Cat_n)$ by invariant theoretic arguments. 
\end{rem}

\begin{lem}
\label{lem:sym1}
For non-negative integers $p, i\geq 0$, 
define a symmetric polynomial $F_{p,i}(x_1, x_2)$ of 
two variables by 
$F_{p,i}:=\frac{x_1^{p+1}-x_2^{p+1}}{x_1-x_2}
\cdot (x_1-x_2)^{2i}$. Then 
$$
\bC[x_1, x_2]^{\frakS_2}=
\bigoplus_{p, i\geq 0}\bC\cdot F_{p,i}. 
$$
\end{lem}
\begin{proof}
Let $F(x_1, x_2)$ be a homogeneous symmetric polynomial. 
We will prove $F$ is expressed a linear combination 
$\{F_{p,i}\}$ by induction on $\deg F$. 
If $\deg F=1$, then $F=x_1+x_2=F_{1,0}$. 
Consider the case $\deg F\geq 2$. If $F(x, x)=0$, 
then we have 
$F=(x_1-x_2)^2\cdot G(x_1, x_2)$ with 
$G$ symmetric. Thus by inductive hypothesis, 
$G$ is a linear combination of $\{F_{p,i}\}$, 
and so is $F$. Suppose $F(x, x)=ax^n\neq 0$. 
Consider $G(x_1, x_2)=F-\frac{a}{n+1}F_{n,0}$. 
Then $G(x, x)=0$, and is reduced to the  previous case. 
Thus $\bC[x_1, x_2]^{\frakS_2}$ is spanned by 
$\{F_{p,i}\}$. By computing the Hilbert series, 
we have (note that $\deg F_{p,i}=p+2i$) 
$$
H(\bC[x_1, x_2]^{\frakS_2}, t)=
\frac{1}{(1-t)(1-t^2)}\leq
H(\sum\bC\cdot F_{p,i}, t)\leq
\frac{1}{(1-t)(1-t^2)}. 
$$
We conclude that $\{F_{p,i}\}$ forms a basis. 
\end{proof}

Let $S=\bC[x_1, \dots, x_n]$. 
%Considering permutations of $\{1, 2\}$ and 
%$\{3, 4, \dots, n\}$, we may 
We consider a subgroup 
$\frakS_2\times\frakS_{n-2}\subset\frakS_n$ 
which acts $\{x_1, x_2\}$ and 
$\{x_3, \dots, x_n\}$ respectively. 
\begin{lem}
\label{lem:sym2}
$$
S^{\frakS_2\times\frakS_{n-2}}=
\bC[x_1, x_2]^{\frakS_2}
\cdot S^{\frakS_n}. 
$$
\end{lem}
\begin{proof}
The inclusion $\supseteq$ is clear. For the 
reverse inclusion, we first note that 
$$
S^{\frakS_2\times\frakS_{n-2}}=
\bC[x_1, x_2]^{\frakS_2}\otimes
\bC[x_3, \dots, x_n]^{\frakS_{n-2}}. 
$$
As is well known, $S^{\frakS_n}$ is 
generated by $x_1^k+\dots+x_n^k$ ($k\geq 0$) as a 
$\bC$-algebra. Thus 
$$
x_3^k+\dots+x_n^k=
(x_1^k+\dots+x_n^k)-(x_1^k+x_2^k)\in
\bC[x_1, x_2]^{\frakS_2}
\cdot S^{\frakS_n}. 
$$
\end{proof}
Combining Lemma \ref{lem:sym1} and \ref{lem:sym2}, 
we have. 
\begin{lem}
\label{lem:sym3}
$$
S^{\frakS_2\times\frakS_{n-2}}=
\sum_{p, i\geq 0}
S^{\frakS_n}\cdot F_{p,i}. 
$$
\end{lem}

Now we prove Theorem \ref{thm:cat}. 
Denote $H_0=\{z=0\}$ and $\alpha_{ij}=x_i-x_j$. 
The restriction $(\A^{H_0}, m^{H_0})$ is equal to 
$(A_{n-1}, 3)$. From the result of \S\ref{sec:coxeter}, 
$(\A^{H_0}, m^{H_0})$ is free with exponents 
$(0, n+1, \dots, 2n-1)$. 
We can choose basis 
$\eta_1, \dots, \eta_n\in D(\A^{H_0}, m^{H_0})^{\frakS_n}$ 
from $\frakS_n$-invariant vector fields. 
By definition, $\eta_i(x_1-x_2)=(x_1-x_2)^3G_i$, and 
$G_i\in S^{\frakS_2\times\frakS_{n-2}}$. By 
Lemma \ref{lem:sym3}, there exist symmetric 
polynomials $B_i^{p,r}\in S^{\frakS_n}$ such that 
\begin{equation}
\label{eq:key}
G_i=\sum_{p, r\geq 0}B_i^{p,r}\cdot F_{p,r}. 
\end{equation}
With $B_i^{p,r}$ and $\frakS_n$-invariant 
vector field $\delta_p=\sum_{i=1}^nx_i^p\partial_i$, 
let us define 
\begin{equation}
\widetilde{\eta}_i:=
\eta_i-\sum_{p,r\geq 0}z^{2r+2}B_i^{p,r}\delta_{p+1}. 
\end{equation}
Then 
$\widetilde{\eta}_1, \dots, 
\widetilde{\eta}_n$ and 
$\theta_E=\sum_{i=1}^nx_i\partial_i+z\partial_z$ 
form a basis of $D(\Cat_n)$. Indeed, 
they are linearly independent over 
$\bC[x_1, \dots, x_n, z]$ (since $\eta_i$'s are 
independent), and 
\begin{eqnarray*}
\widetilde{\eta}_i(\alpha_{12}\pm z)&=&
\widetilde{\eta}_i \alpha_{12}\\
&=&
\eta(x_1-x_2)-
\sum_{p,r\geq 0}z^{2r+2}B_i^{p,r}(x_1^{p+1}-x_2^{p+1})\\
&=&
(x_1-x_2)^3 G_i-
\sum_{p,r\geq 0}z^{2r+2}B_i^{p,r}(x_1^{p+1}-x_2^{p+1}). 
\end{eqnarray*}
The last polynomial is divisible by $z^2-\alpha_{12}^2$. 
Indeed, put $z=\pm\alpha_{12}$ in the last formula, 
we have from (\ref{eq:key})
\begin{eqnarray*}
&=&
(x_1-x_2)^3 G_i-
\sum_{p,r\geq 0}(x_1-x_2)^{2r+2}B_i^{p,r}(x_1^{p+1}-x_2^{p+1})\\
&=&
(x_1-x_2)^3\left(
G_i-
\sum_{p,r\geq 0}(x_1-x_2)^{2r}B_i^{p,r}
\frac{x_1^{p+1}-x_2^{p+1}}{x_1-x_2}
\right)\\
&=&
(x_1-x_2)^3\left(
G_i-
\sum_{p,r\geq 0}B_i^{p,r}F_{p,r}
\right)\\
&=&0. 
\end{eqnarray*}
Since $\widetilde{\eta}_i$ is $\frakS_n$-invariant, 
$\widetilde{\eta}_i\alpha_{jk}$ is divisible 
by $\alpha_{jk}(\alpha_{jk}^2-z^2)$ for 
any $j, k$. This completes the proof.

\section{Concluding remarks and open problems}
\label{sec:pbm}

One of the central problems in the theory of hyperplane 
arrangements is to decide to what extent the 
structure of an arrangement 
is determined by the combinatorics of the arrangement. 
\begin{pbm}
Let $\A_1$ and $\A_2$ be central arrangements in 
$\bK^\ell$. Assume that $L(\A_1)\simeq L(\A_2)$. 
Does the freeness of $\A_1$ imply the freeness of $\A_2$? 
\end{pbm}
The above (by Terao \cite{ter-conj}) is a long 
standing problem, even for the case $\ell=3$, 
since the beginning of this area. 
Note that several variants of this problem are 
known to have counter examples:
\begin{itemize}
\item $\A_1$ and $\A_2$ are multiarrangements 
(\cite{zie-multi}), 
\item $\A_1$ and $\A_2$ are defined over different fields 
(\cite{zie-mat}), 
\item $\A_1$ and $\A_2$ are line-conic arrangements 
(\cite{sch-toh}). 
\end{itemize}

There are several characterizations for $3$-arrangements to 
be free via Ziegler's restriction map 
(\S\ref{subsec:char1} and 
see \cite{abe-exp} for recent developments). 
However the author does not know the answer to the following. 
\begin{pbm}
Does the converse to Theorem \ref{thm:curv} hold? 
\end{pbm}

As in \S\ref{subsec:adj} the modules of 
multiderivations naturally appear in 
the study of relative de Rham cohomology groups for 
the adjoint quotient map $\chi:\frakg\rightarrow\frakh/W$. 
However the idea of the proof 
of Theorem \ref{thm:adjqt} works only for the type $ADE$ and 
the cohomology of degree $2$. 
\begin{pbm}
Study the structure of the relative de Rham cohomology 
group $H^k(\Omega_\chi)$ as $S^W$-module. 
\end{pbm}
\noindent
The structures of $H^2(\Omega_\chi)$ for non simply laced 
cases are expected to be related with the module of derivations 
with a non constant multiplicity. Similarly higher degree 
cases are expected to be related with the module of 
higher derivations 
$D^k(\A, m)$ defined in \cite{atw-char}. 

Postnikov and Stanley \cite{ps-def} observed 
curious properties of the characteristic polynomials 
that for some truncated affine Weyl arrangements, 
the all roots of the characteristic polynomial 
have the same real part (``Riemann hypothesis''). 
By 
Solomon-Terao's formula (Theorem \ref{thm:stf}), the 
characteristic polynomial is determined by the 
module $D(\A)$. It would be natural to expect curious 
behaviours of the characteristic polynomial reflect 
the structure of $D(\A)$. 
The following question posed by Athanasiadis 
(\cite[Question 6.5]{ath-survey}) is still challenging. 
\begin{pbm}
Are there any natural algebraic structures of $D(\A)$ 
which cause Riemann hypothesis? 
\end{pbm}
For instance, is the ``functional equation'' 
(\cite[(9.12)]{ps-def}) deduced from the self duality (up to 
degree shift) 
 $D_0(\A)\simeq D_0(\A)[-d]^\lor$ of certain module?

%%% koko

%%%%%%%%%%%%%%%%%%%%%%%%%%%%%%%%%
% References
%%%%%%%%%%%%%%%%%%%%%%%%%%%%%%%%%

\end{document}